\documentclass[a4paper,12pt]{article}
\usepackage[english]{babel}
\usepackage[T2A]{fontenc}
\usepackage[cp1251]{inputenc}
\usepackage[english]{babel}
\usepackage{amsthm}
\usepackage[tbtags]{amsmath}
\usepackage{amsfonts,amssymb}
\usepackage[usenames]{color}
\usepackage[normalem]{ulem}
\begin{document}

\newtheorem{theorem}{Theorem}
\newtheorem{lemma}{Lemma}
\newtheorem{proposition}{Proposition}
\newtheorem{Cor}{Corollary}

\begin{center}
\textbf{Finite Non-Commutative Centrally Essential Rings of Minimal Order}
\end{center}

\begin{center}
Adel Abyzov\footnote{Kazan (Volga Region) Federal University, Kazan, Russia; email: adel.abyzov@kpfu.ru},
Oleg Lyubimtsev\footnote{Nizhny Novgorod State University, Department of Algebra, Geometry and Discrete Mathematics, scientific and educational mathematical center "Mathematics of technologies of the future", Nizhny Novgorod, Russia; email: oleg\_lyubimcev@mail.ru, https://orcid.org/0000-0002-8773-1943.},
Askar Tuganbaev\footnote{National Research University MPEI, Department of Higher Mathematics, Moscow, Russia; Lomonosov Moscow State University, Department of Higher Algebra, Moscow, Russia; tuganbaev@gmail.com, https://orcid.org/0000-0001-9582-3314.}

\end{center}

\textbf{Abstract.} A ring $R$ with non-zero identity element is said to be centrally essential if for any non-zero element $a$ of $R$, there exist non-zero central elements $x$ and $y$ such that $ax = y$. In the paper, we prove that the minimal order of a finite non-commutative centrally essential ring is equal to $64$. A partial classification of such rings is presented, and the existence of non-isomorphic non-commutative centrally essential rings of minimal order is proved.

\textbf{Keywords:} centrally essential ring, finite local ring.

The work of Adel Abyzov is supported by the Russian Science Foundation (grant 23-21-10086) and performed under the development program of Volga Region Mathematical Center (agreement 075-02-2025-1725/1). 
The work of Oleg Lyubimtsev is supported by Ministry of Education and Science of the Russian Federation, project H-044-0\_2026-2028. The study of Askar Tuganbaev is supported by Russian Science Foundation, project 26-11-00001.

\textbf{MSC2020 subject classification 16D25, 16R99}

\section{Introduction}

We consider only associative rings with a non-zero identity element and unitary modules. A ring $R$ is said to be \textbf{centrally essential} if for any non-zero element $a\in R$, there exist non-zero central elements $x,y\in R$ such that $ax = y$.\footnote{It is clear that the ring $R$ with center $Z$ is centrally essential if and only if the module $R_{Z}$ is an essential extension of the module $Z_{Z}$.} The class of centrally essential rings is a proper extension of the class of commutative rings. In \cite[Proposition 3.3]{MT18}, it is proved that semiprime centrally essential rings are commutative. In \cite[ Example 5.2]{LT25IJ} a non-commutative centrally essential ring of order $128$ is constructed. In this regard, a natural question arises concerning the minimum possible order of a non-commutative centrally essential ring. We denote by $Z(R)$ and $J(R)$ the center and the Jacobson radical of the ring $R$, respectively. It is well known that the Jacobson radical of a right Artinian ring is nilpotent. The smallest $n$ such that $J(R)^n = 0$, is called the \textbf{nilpotence index} of $J(R)$. We write $\mathrm{Ann}(S) = \{r\in R \,|\, rS = 0\}$ for the left annihilator of the set $S$ in the ring $R$. A ring $R$ is said to be \textbf{local} if $R/J(R)$ is a division ring.

The main results of the paper are Theorems 1.1, 1.2 and 1.3.

\textbf{1.1. Theorem.} The smallest finite non-commutative centrally essential ring has $64$ elements. 

\textbf{1.2. Theorem.} Let $R$ be a finite non-commutative centrally essential ring of order $64$. Then:
\begin{enumerate}
 \item The ring $R$ is local.
 \item The factor ring $R/J(R)$ is isomorphic to the field $\mathbb{F}_2$ and the radical $J(R)$ has order $32$.
 \item The center $Z(R)$ has order $16$ and the intersection $Z_0 = Z(R) \cap J(R)$ has order $8$. 
\end{enumerate}

\textbf{1.3. Theorem.} Within the class of non-commutative centrally essential algebras of order $64$ over the field $\mathbb{F}_2$ satisfying the condition 
$J(R)^2 \subseteq \mathrm{Ann}(J(R)) \cap Z(R)$, there exist exactly two families distinguished by an algebraic invariant -- the dimension of the annihilator $\dim_{\mathbb{F}_2} W_0 \in \{2, 1\}$, where $W_0 = \mathrm{Ann}(J(R)) \cap Z(R)$. These families are described as follows.

\textbf{Family 1 ($\dim_{\mathbb{F}_2} W_0 = 2$)}

The algebra has a basis $\{1, a, b, c, s_1, s_2\}$, with $J(R) = \langle a, b, c, s_1, s_2 \rangle_{\mathbb{F}_2}$, $Z(R) = \langle 1, c, s_1, s_2 \rangle_{\mathbb{F}_2}$, and $W_0 = \langle s_1, s_2 \rangle_{\mathbb{F}_2}$. The multiplication of the basis elements of the radical is defined by the rules: $ac = ca = s_1$; $bc = cb = s_2$; $ab = s_1$, $ba = 0$; $c^2 = 0$. The elements $s_1$ and $s_2$ annihilate $J(R)$. The squares $a^2 = \alpha$ and $b^2 = \beta$ are   parametrized by elements $\alpha, \beta\in W_0$.

\textbf{Family 2 ($\dim_{\mathbb{F}_2} W_0 = 1$)} 

The algebra has a basis $\{1, a, b, c_1, c_2, w\}$, 
$J(R) = \langle a, b, c_1, c_2, w \rangle_{\mathbb{F}_2}$, 
$Z(R) = \langle 1, c_1, c_2, w \rangle_{\mathbb{F}_2}$ and 
$W_0 = \langle w \rangle_{\mathbb{F}_2}$.
The multiplication of basis elements of the radical is defined by the rules:
$ac_1 = c_1a = w$, $bc_1 = c_1b = 0$;
$ac_2 = c_2a = 0$, $bc_2 = c_2b = w$;
$ab = w$, $ba = 0$; $c_1^2 = c_2^2 = c_1 c_2 = c_2 c_1 = 0$.
The element $w$ annihilates the radical $J(R)$. The squares $a^2 = \alpha w$ and $b^2 = \beta w$ are parametrized by scalars $\alpha, \beta \in \mathbb{F}_2$.

\section{Proof of Theorem 1.1}

Theorem 1.1 follows from Lemmas 2.1, 2.2 and 2.3.

\textbf{2.1. Lemma.} There do not exist finite non-commutative centrally essential rings of order $\le 32$. 

\textbf{Proof.} Let $R$ be a finite non-commutative centrally essential ring of minimal order. Then $R$ is a local ring. Indeed, in any centrally essential ring, all idempotents are central; see \cite[Lemma 2.3]{MT18}. Lifting idempotents modulo a nilpotent ideal allows a unique decomposition of a finite ring $R$ into a direct product of finite local rings. It is clear that a direct product of such rings is centrally essential if and only if each of its  factors is a centrally essential ring. If there were more than one factor, then a non-commutative factor would have strictly smaller order, which contradicts the minimality of $R$. 

Let $|R| = p^n$, $Z = Z(R)$, $J = J(R)$ and the factor ring $R/J \cong \mathbb{F}_q$, where $q = p^f$. We set $Z_0 = Z \cap J$. Let us obtain an upper bound for the cardinality of the set $Z_0$. Since the ring $R$ is not commutative, the index of its center as an additive subgroup, $[R : Z]$, is at least $p^2$. Indeed if $[R : Z] = p$, then the additive factor group $R/Z$ is a cyclic group generated by some residue class $a+Z$. Then any element of the ring is of the form $z + ka$, where $z \in Z$, 
$k \in \mathbb{Z}$ and $[z_1 + k_1a, z_2+k_2a] = 0$, i.e. the ring $R$ is commutative. Consequently, $|Z| \le |R|/p^2 = p^{n-2}$. Since $1 \in Z$ and $1 \notin J$, we have that 
$Z_0 = Z \cap J$ is a proper subgroup in $Z$. Since $Z$ is a finite $p$-group, the index of any of its proper subgroup is not less than $p$. Therefore, $|Z| \ge p|Z_0|$. Combining the inequalities, we obtain an upper bound for any prime integer $p$:
$$
|Z_0| \le \frac{|Z|}{p} \le \frac{p^{n-2}}{p} = p^{n-3}.
$$ 
Now we obtain a lower bound on the cardinality of the set $Z_0$. First, note that if $J = 0$, then the ring $R$ is semiprime. However, semiprime centrally essential rings are commutative; see \cite[Proposition 3.3]{MT18}. Consequently, $J$ is a nilpotent ideal with index of nilpotency $k$, where $k \geq 2$. Since $J$ is a nilpotent ideal with nilpotence index $k$, it follows from \cite[Lemma 2.4]{AS01} that $J^{k-1}\cap Z(R)$ is a central ideal. In addition, since the ring $R$ is centrally essential, we have $J^{k-1}\cap Z(R)\neq 0$. Therefore, the annihilator $\mathrm{Ann}(J)$ contains non-zero central element $z_1 \in W_0 = \mathrm{Ann}(J) \cap Z_0$. Since the ring $R$ is not commutative, it follows from \cite{H54} that there exists an element 
$r \in J \setminus Z_0$. Then there exists an element $x \in Z \setminus \{0\}$ such that 
$rx \in Z \setminus \{0\}$. If $x \notin J$, then $x$ is invertible in $R$, whence 
$r \in Z x^{-1} = Z$ and $r \in Z_0$. This is a contradiction. Therefore, 
$x \in Z_0 \setminus \{0\}$. Since $r z_1 = 0$ and $rx \neq 0$, we have that $x \neq z_1$. Consequently, $Z_0$ contains at least three distinct elements $0$, $z_1$ and $x$. We obtain that 
$|Z_0| \ge 3$. We also remark that for $p = 2$, we have that $|Z_0| \ge 4$, since the order is a power of two. Consequently, $3 \le |Z_0| \le p^{n-3}$. We verify possible orders $|R| \le 32$. If $p \ge 3$, it follows from $|R| \le 32$ that $n \le 3$. Then $|Z_0| \le p^{3-3} = 1$; this contradicts the fact that $|Z_0| \ge 3$. Therefore, the characteristic of the residue field cannot be odd. 
 
Let $p = 2$. We have $|Z_0| \le 2^{n-3}$. Since $|Z_0| \ge 4$, we have that $n \ge 5$. By considering $|R| \le 32$, we obtain $n = 5$, that is $|R| = 32$. In addition, $|Z_0| \le 2^{5-3} = 4$, whence $|Z_0| = 4$. We consider possible variants for the degree $f$ (where $q = 2^f$). In a finite local ring of order $p^n$ with residue field $\mathbb{F}_{p^f}$, the degree $f$ must be a divisor of $n$. Indeed, since the characteristic of the field $\mathbb{F}_{p^f}$ is $p$, we have 
$p \cdot 1 \in J$, which implies $p J^k \subseteq J^{k+1}$ for all $k \ge 0$. Therefore, the filtration quotients $J^k/J^{k+1}$ are vector spaces over the field $R/J \cong \mathbb{F}_{p^f}$. It follows that $f$ is a divisor of $n$. 
Since $n = 5$ is a prime integer, cases $f = 1$ and $f = 5$ are only possible. If $f = 5$, then 
$|J| = 32/32 = 1$ and the ring $R$ is commutative.  This is a contradiction. Therefore, the unique possible variant is $f=1$ (in addition $|J|=16$). 

The additive group, $Z_0$ consists of $4$ elements: $\{0, z_1, x, y\}$, where $z_1\in W_0$.  Since $W_0 \subsetneq Z_0$, we have $W_0 = \{0, z_1\}$. Then, for all $r \in J \setminus Z_0$, we have
$0 \neq ry = r(x + z_1) = rx \in Z_0$ since  the ring $R$ is centrally essential.

We consider  a homomorphism of additive groups $L_x\colon J \to Z_0$ such that $r \mapsto rx$. It follows from the above and the inclusion $Z^2_0 \subseteq Z_0$ that $L_x$ is well-defined and $\ker(L_x) \subseteq Z_0$. Since $|\text{Im}(L_x)| \le |Z_0| = 4$, we have $|J|/|\ker(L_x)| \le 4$, whence $|\ker(L_x)| \ge 16/4 = 4$. Thus, $\ker(L_x) = Z_0$, and $L_x$ induces a bijection $J/Z_0 \to Z_0$. There exists $u \in J$ such that $ux = x$, implying $(1-u)x = 0$. Since $u \in J$, the element $(1 - u)$ is invertible. Hence, $x = 0$. This is a contradiction. Thus, all possible cases for $|R| \le 32$ are ruled out, and the lemma is proved.~$\square$

\textbf{2.2. Lemma.}
There are no finite non-commutative centrally essential rings with orders strictly greater than 32 and strictly less than 64. 

\textbf{Proof.} Let $R$ be a finite non-commutative centrally essential ring, $32 < |R| < 64$. A ring $R$ is decomposed into a direct product of finite local centrally essential rings:
$$
R = R_1 \times R_2 \times \dots \times R_k.
$$
Let $R_1$ be a non-commutative factor. If $k \ge 2$, then the order of the factor $|R_1|$ is a proper divisor of order $|R|$. The greatest proper divisor of any positive integer $\le 63$ does not exceed $[63/2] = 31$. Therefore, $|R_1| \le 31$. By Lemma 2.1, such a ring does not exist. If $k = 1$, then the ring $R$ is local and its order must be a power of prime integer.  Among the powers of prime numbers in the interval $(32, 64)$, there are only prime numbers and the square of a prime number, $49 = 7^2$. Every ring of prime order or of order $p^2$ is always commutative; see \cite{E68}. This is a contradiction.~$\square$

\textbf{2.3. Lemma.} There exists a non-commutative centrally essential ring consisting of $64$ elements. 

\textbf{Proof.} We construct a $6$-dimensional algebra $R$ over the field of two elements 
$\mathbb{F}_2$. The order of such an algebra is equal to $|R| = 2^6 = 64$. We define the following basis of the algebra: 
$\{1, a, b, c, s_1, s_2\}$. The unity $1$ is a multiplicative identity element. The Jacobson radical
$J $ is equal to $\langle a, b, c, s_1, s_2 \rangle_{\mathbb{F}_2}$. For basis elements of the radical, we set
$$
a^2 = b^2 = c^2 = 0, ab = s_1, ba = 0, ac = ca = s_1, bc = cb = s_2,
$$
$$
s_1 x = x s_1 = 0, s_2 x = x s_2 = 0 \quad \mbox{for all} \, x\in J.
$$
Since the product of any two elements from $J$ belongs to the subspace 
$S = \langle s_1, s_2 \rangle_{\mathbb{F}_2}$ and $S\subseteq \mathrm{Ann}(J)$, we have that $J^3 = 0$ and the identity $(xy)z = x(yz)$ is trivially true. Next, since $ab = s_1\neq 0$ and $ba = 0$, we have that the ring $R$ is not commutative.

Let us find the center $Z(R)$ of the ring $R$. Let $x \in R$, then $x = k_0 1 + k_a a + k_b b + k_c c + k_1 s_1 + k_2 s_2$, where $k_i \in \mathbb{F}_2$. We have:
$$
[x, a] = xa - ax = k_b[b,a] + k_c[c,a] = k_b(0 - s_1) + k_c(s_1 - s_1) = k_b s_1,
$$
whence $k_b = 0$. Next,
$$
[x, b] = xb - bx = k_a[a,b] + k_c[c,b] = k_a s_1 + k_c(s_2 - s_2) = k_a s_1,
$$
whence $k_a = 0$. Since $[x, c] = 0$, we have that $Z(R) = \langle 1, c, s_1, s_2\rangle_{\mathbb{F}_2}$.

Every non-central element $r \in R \setminus Z(R)$ is of the form $r = k_a a + k_b b + z$, where 
$z \in Z(R)$ and $(k_a, k_b) \neq (0,0)$. Then
$$
rc = (k_a a + k_b b + k_0 1 + k_c c + k_1 s_1 + k_2 s_2)c = k_0 c + k_a s_1 + k_b s_2.
$$
Since the basis vectors $c, s_1, s_2$ are linearly independent and at least one of the coefficients $k_a, k_b$ is non-zero, it follows that $rc \neq 0$ and $rc \in Z(R)$. Therefore, the algebra $R$ is centrally essential. ~$\square$

\section{Proof of Theorem 1.2}

\indent
{\bf 1.} Since all idempotents in a centrally essential ring are central, a finite centrally essential ring decomposes uniquely into a direct product of local centrally essential rings. If the ring $R$ has a non-trivial decomposition, then the order of the non-commutative factor does not exceed $64/2 = 32$. However, by Lemma 2.1, there are no non-commutative centrally essential rings of order $\le 32$. Consequently, the ring $R$ is local. 

{\bf 2.} Let $R/J(R) \cong \mathbb{F}_q$.  Since $|R| = 64 = 2^6$, we have $q = 2^f$. Thus, $f$ must be a divisor of $6$, i.e., $f \in \{1, 2, 3, 6\}$. If $q = 64$ ($f = 6$), then $|J(R)| = 1$. Consequently, the ring $R \cong \mathbb{F}_{64}$ is commutative, which contradicts the theorem's hypothesis. We prove that $q = 2$ ($f = 1$). We assume that $q \ge 4$. It follows from $[R:Z(R)] \ge 4$ that $|Z(R)| \le 16$. Since $[Z(R):Z_0] \ge 2$, we have that $|Z_0| \le 8$. The set $W_0 = \mathrm{Ann}(J) \cap Z_0$ is non-zero (see the proof of Lemma 2.1), consequently $|W_0| \ge 2$. Since the ring $R$ is centrally essential, we have that $W_0 \subsetneq Z_0$, whence $|Z_0| \ge 4$. 

Let $q = 8$ ($f = 3$). Then $|J(R)| = 8$. From the estimate obtained earlier we have 
$4 \le |Z_0| \le 8$. If $|Z_0| = 8$, then $Z_0 = J(R) \subseteq Z(R)$ and the ring $R$ is commutative. If $|Z_0| = 4$, then $|J(R)/Z_0| = 2$. Therefore, the radical 
$J(R) = \langle Z_0, j_0 \rangle$ is commutative. The canonical epimorphism of the groups of units $R^* \to (R/J(R))^* \cong \mathbb{F}_8^*$ has a kernel $1 + J(R)$ of order $8$. Consequently, there exists an element $\omega \in R^*$ of multiplicative order $7$ such that its coset $\omega + J(R)$ generates the field $R/J(R)$. The action by conjugation by $\omega$ on the two-element quotient group $J(R)/Z_0$ is trivial; therefore, $\omega j_0 \omega^{-1} = j_0 + z$ for some $z \in Z_0$. Multiplying on the right by $\omega$, we obtain $\omega j_0 = j_0 \omega + z \omega$. Since 
$z \in Z(R)$, applying this relation repeatedly yields $\omega^k j_0 = j_0 \omega^k + k z\omega^k$. For $k=7$, we have $\omega^7 = 1$, which implies $j_0 = j_0 + 7z$. Because $Z_0$ is a finite $2$-group, multiplication by the odd number $7$ is an automorphism of the group $Z_0$. It follows that $z = 0$, and thus $\omega$ commutes with $j_0$. Since $\omega$ and $J(R)$ generate the entire ring $R$, and $J(R)$ is commutative, the ring $R$ is commutative.
This is a contradiction.
 
Suppose $q = 4$ ($f = 2$). Then $|J(R)| = 16$. Let us consider the possible orders 
$|Z_0| \in \{4, 8\}$. Suppose $|Z_0| = 4$. Then $|W_0| = 2$. Choose an element 
$c \in Z_0 \setminus W_0$. Following the same steps as in the proof of Lemma 2.1, one can verify that the homomorphism $L_c: J(R) \to Z_0$ (defined by $r \mapsto rc$) induces a bijection 
$J(R)/Z_0 \to Z_0$, which leads to a contradiction. Finally, assume $|Z_0| = 8$. Then 
$|J(R)/Z_0| = 16/8 = 2$, which implies that $J(R) = \langle Z_0, j_0 \rangle$ is a commutative ideal. Analogous to the case $q=8$, choose an element $\omega \in R^*$ of multiplicative order $3$ such that its image generates $R/J(R) \cong \mathbb{F}_4$. This element acts trivially on $J(R)/Z_0$: $\omega j_0 \omega^{-1} = j_0 + z$ for some $z \in Z_0$. Raising this to the third power yields $j_0 = j_0 + 3z$, which implies $z=0$ and $\omega j_0 = j_0 \omega$. Since $\omega$ and $J(R)$ generate $R$, the ring $R$ is commutative. This is a contradiction. 

Thus, the assumption $q \ge 4$ leads to a contradiction in all cases. Consequently, 
$R/J(R) \cong \mathbb{F}_2$ and $|J(R)| = 64/2 = 32$. 

{\bf 3.} Since the ring $R$ is not commutative, the index of its center as an additive group, $[R : Z(R)]$, is at least 4. Therefore, $|Z(R)| \le 64/4 = 16$. It follows from $1 \in Z(R) \setminus J(R)$ that $Z_0 = Z(R) \cap J(R)$ is a subgroup of index $2$ in $Z(R)$. Consequently, $|Z_0| = |Z(R)| / 2 \le 8$. Since 
$W_0 \subsetneq Z_0$ and $W_0 \neq \{0\}$, we have that $|W_0| \ge 2$. Therefore, $|Z_0| \ge 4$. This implies that $|Z_0| = 4$ or $|Z_0| = 8$; in addition $|Z(R)| = 8$ or $|Z(R)| = 16$, respectively. Suppose that $|Z(R)| = 8$ and $|Z_0| = 4$, $Z_0 = \{0, z_1, x, y\}$, where $z_1\in W_0$ (see the proof of the lemma 2.1). Then the mapping $L_x : J(R) \to Z_0$ ($r \mapsto rx$) is a well-defined homomorphism of additive groups with kernel $\ker(L_x) \subseteq Z_0$. But then we have $|J(R)| = |ker(L_x)| Ў¤ |Im(L_x)| \le |Z_0|\cdot|Z_0| = 4 \cdot 4 = 16$, which contradicts the equality $|J(R)| = 32$.~$\square$

\section{Proof of Theorem 1.3}

We set $J = J(R)$ and $Z_0 = Z(R) \cap J$. Since $J^2 \subseteq W_0$, we have that $J^3 = 0$ and the associativity is proved. We have $\dim J = 5$. The case $\dim Z_0 = 2$ is impossible. Indeed, if 
$\dim Z_0 = 2$, then $\dim W_0 = 1$ (since $0 \subsetneq W_0 \subsetneq Z_0$). Let 
$Z_0 = \langle c, w \rangle$ and $W_0 = \langle w \rangle$. We consider the operator of right multiplication $L_c: J \to W_0$ defined by the rule $r \mapsto rc$. It is well defined, since for any 
$r \in J$, we have that $rc \in J^2 \subseteq W_0$. Since $\dim J = 5$ and $\dim W_0 = 1$, it follows from the rank-nullity theorem $\dim \ker L_c \ge 5 - 1 = 4$. Consequently, $|\ker L_c| \ge 16$. Since 
$|Z_0| = 4$, there exists a non-central element $r \in \ker L_c \setminus Z_0$ and $rc = 0$. 
Since $w \in W_0 \subseteq \mathrm{Ann}(J)$, we have that $rw = 0$. Therefore, $rZ_0 = 0$; this contradicts the fact that the ring $R$ is centrally essential. Consequently, $\dim Z_0 = 3$, whence 
$\dim W_0 \in \{1, 2\}$. We consider these cases. 

\textbf{Family 1 ($\dim W_0 = 2$).}
Let $W_0 = \langle s_1, s_2 \rangle$ and $Z_0 = \langle c, s_1, s_2 \rangle$. We consider the operator of right multiplication $L_c: J \to W_0$ ($r \mapsto rc$). Since $\dim J = 5$ and 
$\dim W_0 = 2$, the dimension of the kernel $\dim \ker L_c$ is not less than $3$. In addition, $rZ_0\neq 0$ for the non-central element $r$, otherwise, the ring $R$ is not centrally essential. Then $\ker L_c \subseteq Z_0$. By considering $\dim Z_0 = 3$, we obtain $\ker L_c = Z_0$, which means that $Z_0 c = 0$. From this, $Z^2_0 = 0$ follows, given that $Z_0W_0 = 0$ (in particular, $c^2 = 0$). Consequently, $L_c$ induces a well-defined isomorphism $J/Z_0 \to W_0$. We choose a elements $a, b \in J \setminus Z_0$ whose cosets form a basis of $J/Z_0$ such that $ac = s_1$ and $bc = s_2$. 

Since the algebra $R$ is not commutative, the commutator $[a,b]$ is not equal to $0$. Given that $J^2 \subseteq W_0$, we have $t = [a,b] \in W_0 \setminus \{0\}$. We show that, by suitably changing the bases in $J$, we can achieve the equality $[a,b] = s_1$. If $t = s_1$, the required equality holds. We show that if $t = s_2$, we can satisfy the condition $[a', b'] = s_1'$. Let us make the substitution $a' = b$, $b' = a + b$. Since the characteristic of the ring $R$ is $2$, we have that
$$
[a', b'] = [b, a + b] = [b, a] = -[a, b] = s_2.
$$
Next, we pass to a new basis for the space $W_0$ by setting $s_1' = s_2$ and $s_2' = s_1 + s_2$. In the new basis, we obtain the desired relation $[a', b'] = s_1'$. Multiplication by $c$ is also consistent with the new basis: $a' c = b c = s_2 = s_1'$ and $b' c = (a+b) c = s_1 + s_2 = s_2'$. 
In the case $t = s_1 + s_2$, by making the substitutions $a' = a + b$, $b' = a$ and 
$s_1' = s_1 + s_2$, $s_2' = s_1$, we again obtain $[a', b'] = s_1'$. Thus, we may assume that 
$[a, b] = s_1$.

In addition, without loss of generality, we can set $ba = 0$. Indeed, given that $Z_0^2 = 0$, upon substituting $a' = a + x$ and $b' = b + y$ (where $x, y \in Z_0$), we have $b'a' = ba + xb + ya$. Since the operators $L_a$ and $L_b$ map $Z_0$ onto the one-dimensional subspaces 
$\langle s_1 \rangle$ and $\langle s_2 \rangle$, respectively, we have that the elements of $W_0$ are linear combinations $xb + ya$. Therefore, there exist $x,y \in Z_0$ such that $xb + ya = ba$ and $b'a' = 0$. In addition, $a'c = (a+x)c = ac + xc = s_1 + 0 = s_1$, since $xc = 0$ for $x \in Z_0$. Next, $b'c = (b + y)c = bc + yc = s_2 + 0 = s_2$. Consequently, the substitution is correct. The squares $a^2 = \alpha$ and $b^2 = \beta$ are selected from $W_0$. Distinct pairs $(\alpha, \beta)$ parametrize this family.

Note that under this substitution, the commutator $[a', b'] = [a, b] = s_1$ remains unchanged. Finally, since the squares $a^2 = \alpha$ and $b^2 = \beta$ belong to $W_0$, they do not affect the commutators of the basis elements or the multiplication by the element $c$. Therefore, as in the proof of Lemma 2.3, for any $\alpha, \beta \in W_0$, the equalities $Z(R) = \langle 1, c, s_1, s_2 \rangle_{\mathbb{F}_2}$ and $W_0 = \langle s_1, s_2 \rangle_{\mathbb{F}_2}$ hold; furthermore, for any non-central element $r \in R \setminus Z(R)$, we have $0 \neq rc \in Z(R)$, which means that the algebra $R$ is centrally essential.

\textbf{Family 2 ($\dim W_0 = 1$).}
Let $W_0 = \langle w \rangle$ and $Z_0 = \langle c_1, c_2, w \rangle$. Since $J^2 \subseteq W_0$ and $JW_0 = 0$, the following linear map is well-defined:
$$
\Phi: J \to \mathrm{Hom}_{\mathbb{F}_2}(Z_0/W_0, W_0), \qquad \Phi(r)(z + W_0) = rz.
$$
If $r \in \ker \Phi$, then $rZ_0 = 0$. The central essentiality of the algebra $R$ implies that no non-central element $r \in J \setminus Z_0$ can annihilate $Z_0$; hence, $\ker \Phi \subseteq Z_0$. Given that $\dim J = 5$ and $\dim \mathrm{Hom}_{\mathbb{F}_2}(Z_0/W_0, W_0) = 2$, the rank-nullity theorem yields $\dim \ker \Phi \ge 5 - 2 = 3$. Since $\dim Z_0 = 3$, it follows that
$\ker \Phi = Z_0$, which implies $Z_0^2 = 0$, that is, $c_1^2 = c_2^2 = c_1c_2 = c_2c_1 = 0$. Furthermore, $\Phi$ induces an isomorphism between the two-dimensional spaces
$\bar{\Phi}: J/Z_0 \to \mathrm{Hom}_{\mathbb{F}_2}(Z_0/W_0, W_0)$. Consequently, one can choose elements $a, b \in J \setminus Z_0$ whose images in $J/Z_0$ form a basis dual to the basis 
$\{c_1 + W_0, c_2 + W_0\}$ of the space $Z_0/W_0$, such that
$$
ac_1 = w, \quad ac_2 = 0, \quad bc_1 = 0, \quad bc_2 = w.
$$
Since the algebra $R$ is non-commutative, the commutator $[a,b] \neq 0$. Given that 
$[a,b] \in J^2 \subseteq W_0 = \langle w \rangle$, we have $ab + ba = [a,b] = w$. Since 
$ab, ba \in \{0, w\}$, there are two possibilities: $ab = w, ba = 0$ or $ab = 0, ba = w$. In the latter case, we can perform the substitution $a' = a + c_2$. Since $Z_0^2 = 0$, the relations $a'c_1 = w$ and $a'c_2 = 0$ are preserved, while $a'b = ab + c_2b = w$ and $ba' = ba + bc_2 = 0$. Thus, without loss of generality, we may assume that $ab = w$ and $ba = 0$. The squares 
$a^2 = \alpha w$ and $b^2 = \beta w$ are parameterized by scalars $\alpha, \beta \in \mathbb{F}_2$. 

Conversely, for any pair $(\alpha, \beta) \in \mathbb{F}_2^2$, the algebra $R$ defined by these rules has a center $Z(R) = \langle 1, c_1, c_2, w \rangle_{\mathbb{F}_2}$ and 
$W_0 = \langle w \rangle_{\mathbb{F}_2}$. For an arbitrary non-central element 
$r = k_0 1 + k_a a + k_b b + z \in R \setminus Z(R)$ (where $z \in Z_0$ and 
$(k_a, k_b) \neq (0,0)$), we have $rc_1 = k_0 c_1 + k_a w$ and $rc_2 = k_0 c_2 + k_b w$. Since 
$(k_a, k_b) \neq (0,0)$, at least one of the products $rc_1$ or $rc_2$ is non-zero and belongs to $Z_0 \subseteq Z(R)$, which implies that $R$ is centrally essential.~$\square$

\textbf{4.1. Remark.} 
Distinct pairs $(\alpha, \beta)$ define algebras within Family 1; their classification (up to isomorphism) requires the construction of precise invariants, leaving this subclass open for further study. For \textbf{Family 2} ($\alpha, \beta \in \mathbb{F}_2$), a complete classification up to isomorphism is obtained quite easily. Specifically, for any element $x = x_1 a + x_2 b + z \in J(R)$ (where $x_1, x_2 \in \mathbb{F}_2$ and $z \in Z_0$), given the centrality of $Z_0$, the identity $Z_0^2 = 0$, and the characteristic $2$, we have:
$$
x^2 = x_1^2 a^2 + x_2^2 b^2 + x_1 x_2 (ab + ba) = (x_1^2 \alpha + x_2^2 \beta + x_1 x_2) w = (x_1 \alpha + x_2 \beta + x_1 x_2) w.
$$ 
Furthermore, the quadratic function $Q(x_1, x_2) = x_1 \alpha + x_2 \beta + x_1 x_2$ over the field $\mathbb{F}_2$ vanishes at three pairs $(x_1, x_2) \in \mathbb{F}_2^2$ when $(\alpha, \beta) \in \{(0,0), (1,0), (0,1)\}$ (in which case $J(R)$ contains exactly $3 \cdot |Z_0| = 24$ elements with a zero square), and only at the zero pair $(0,0)$ when $(\alpha, \beta) = (1,1)$ (in which case $J(R)$ contains exactly $1 \cdot |Z_0| = 8$ elements with a zero square). Moreover, the algebras with parameters $(1,0)$ and $(0,1)$ are isomorphic to the algebra with parameters $(0,0)$. If $\alpha = 1, \beta = 0$ ($a^2 = w, b^2 = 0$), then the substitution
$$
a' = b + c_1, \quad b' = a + b, \quad c_1' = c_1 + c_2, \quad c_2' = c_1
$$
preserves all the defining relations of Family 2: 
$a'c_1' = w$, $a'c_2' = 0$, $b'c_1' = 0$, $b'c_2' = w$, $a'b' = w$, $b'a' = 0$, while yielding 
$(a')^2 = b^2 = 0$ and$(b')^2 = a^2 + b^2 + ab + ba = w + 0 + w + 0 = 0$, i.e., 
$(\alpha', \beta') = (0,0)$. Similarly, the case $(0,1)$ reduces to $(0,0)$ via the substitution $a' = a$, $b' = a + b + c_2$, $c_1' = c_1 + c_2$, $c_2' = c_2$. Thus, Family 2 contains exactly two pairwise non-isomorphic algebras, corresponding to the parameter pairs $(\alpha, \beta) = (0,0)$ and $(\alpha, \beta) = (1,1)$.

\textbf{4.2 Remark.} There exists a non-commutative centrally essential ring of order $64$ which is not an algebra over the field $\mathbb{F}_2$ and has characteristic $4$.

Indeed, we consider the additive group 
$R \cong \mathbb{Z}_4 \oplus \mathbb{Z}_2 \oplus \mathbb{Z}_2 \oplus \mathbb{Z}_2 \oplus \mathbb{Z}_2$ of order $64$. We define generators $1$ (it has order 4) and $a, b, c, s$ (they are of order 2). The element $2 \cdot 1 = 2$ has order $2$ and belongs to $J(R)$. Let us find the radical $J = \langle2, a, b, c, s \rangle_{\mathbb{Z}_2}$ of order $32$. We define the multiplication of elements of $J$ by the rules $a^2 = b^2 = c^2 = 0$, $ab = 2, ba = 0$, $ac = ca = 2$, $bc = cb = s$.
We remark that the elements $2$ and $s$ annihilate $J(R)$. Since $J^3 = 0$, we have that the constructed ring is associative. In addition, it follows from $ab = 2 \neq 0 = ba$ that the ring $R$ is not commutative. Since $2 \cdot 1 \neq 0$, the characteristic of the ring $R$ is equal to $4$. Consequently, $R$ is not an algebra over $\mathbb{F}_2$. Let us find the center $Z(R)$. To do this, we calculate commutators of the element $x = k_0 1 + k_a a + k_b b + k_c c + k_s s$, where $k_0 \in \mathbb{Z}_4$, 
$k_a, k_b, k_c, k_s \in \{0,1\}$:
$$
[x, a] = k_b[b,a] + k_c[c,a] = k_b(-2) = 2k_b \quad \mbox{(calculation in} \, \mathbb{Z}_4).
$$
Then $k_b = 0$. Next, $[x, b] = k_a[a,b] = 2k_a$. Therefore, $k_a = 0$. Consequently, the center 
$Z(R)$ coincides with the set $\langle 1, c, s \rangle$. The intersection $Z_0$ of $Z(R)$ with the radical $J$ is equal to $ \langle 2, c, s \rangle$; it has order $8$. 

For any non-central element $r \in J \setminus Z_0$, we have that $r = k_a a + k_b b + z$, where 
$z \in Z_0$ and $(k_a, k_b) \neq (0,0)$. Then $rc = k_a(ac) + k_b(bc) + zc = k_a(2) + k_b(s)$, where
$c \in Z(R) \setminus \{0\}$. Since the elements $2$ and $s$ are linearly independent over 
$\mathbb{Z}_2$, their linear combination is equal to zero if and only if $k_a = 0$ and $k_b = 0$. Since the element $r$ is not central, we have that $(k_a, k_b) \neq (0,0)$. Consequently, we have $rc \neq 0$. In addition, $rc \in Z_0 \subseteq Z(R)$. Next, for $r = 1 + j \notin J$, where 
$k_0\in \{1, 3\}$, the multiplication by $c$ implies $rc = c + jc$. Since 
$c \notin W_0 = \langle 2, s \rangle_{\mathbb{Z}_2}$ and $jc \in J^2 \subseteq W_0$, their sum is not equal to zero. Therefore, $rc \in Z_0 \setminus \{0\}$. The centrally essential condition is satisfied in all cases, and the ring with the structure specified in the remark exists.

\end{document}